\documentclass[11pt]{article}

\usepackage[T1]{fontenc}
\usepackage{lmodern}
\usepackage[margin=1in]{geometry}
\usepackage{microtype}
\usepackage{amsmath,amssymb,amsthm,mathtools}
\usepackage{xcolor}
\usepackage[hidelinks]{hyperref}
\usepackage{comment}
\usepackage{soul}
\usepackage{graphicx}
\setstcolor{red} 
\setul{}{1pt} 

\hypersetup{
  pdftitle={Memoryless BFGS on a nonsmooth convex function: two new results},
  pdfauthor={Azam Asl},
  pdfsubject={An addendum to the analysis of memoryless BFGS on a nonsmooth convex function}
}
\newtheorem{theorem}{Theorem}[section]

\theoremstyle{remark}

\numberwithin{equation}{section}
\newcommand{\R}{\mathbb{R}}
\DeclareMathAlphabet{\mathbbold}{U}{bbold}{m}{n}
\newcommand{\one}{\mathbbold{1}}
\newcommand{\abs}[1]{\lvert#1\rvert}
\newcommand{\norm}[1]{\lVert#1\rVert}
\DeclareMathOperator{\sgn}{sgn}
\allowdisplaybreaks[1]

\title{Memoryless BFGS on a nonsmooth convex function:\\
two new results}
\author{Azam Asl \and Michael L. Overton}
\date{September 2026}

\begin{document}
\maketitle

\begin{abstract}
We extend the analysis of Asl and Overton~\cite{AO2021} for memoryless BFGS (also known as L-BFGS-1)
using an Armijo-Wolfe line search and applied to the nonsmooth convex function
$f(x)=a\abs{x^{(1)}}+\sum_{i=2}^n x^{(i)}$, where $n\ge2$ and $a \geq \sqrt{n-1}$.
This function is unbounded below, so the method is successful if it generates function values diverging to $-\infty$,
and fails otherwise. 
For scaled L-BFGS-1, we sharpen the sufficient condition for 
failure in
\cite{AO2021} from $a\ge2\sqrt{n-1}$ to $a\ge\sqrt{3(n-1)}$.

For the unscaled version, we prove success for every $a$: either the method generates a direction along which no Armijo-Wolfe step exists, certifying unboundedness below, or its function values diverge to $-\infty$.

Both results are consistent with the numerical experiments in \cite{AO2021}.
Both proofs given here were derived with help from GPT 6 in discussions directed by the first author and corrected and rephrased afterward.
\end{abstract}


\section{Introduction}
\label{sec:introduction}

Asl and Overton~\cite{AO2021} studied the effect of scaling on memoryless
BFGS, also called L-BFGS-1, using an Armijo-Wolfe line search, applied to the nonsmooth function
\begin{equation}
 f(x)=a\abs{x^{(1)}}+\sum_{i=2}^n x^{(i)},\qquad
  n\ge2,\quad  r:=n-1,\quad a \ge \sqrt{r}.
 \label{eq:model}
\end{equation} 
This function is convex and unbounded below. It is differentiable precisely
when $x^{(1)}\ne0$, with
\begin{equation}
 g(x):=\nabla f(x)=
 \left[\begin{array}{c}a\,\sgn(x^{(1)})\\ \one\end{array}\right],
 \qquad \one\in\R^r.
 \label{eq:gradient}
\end{equation}
Here 
$\one$ denotes the all-ones vector. 
Note that although $f$ is unbounded below, it is bounded below along the line defined
by the negative gradient direction from any point $x$ with $x^{(1)} \ne0$, since $a \geq \sqrt{r}$. The case $a < \sqrt{r}$
is of no interest since $f$ is unbounded below along such negative gradient directions.
Throughout the rest of this article we continue to use the same definitions and notation as in \cite{AO2021} and introduce a notation only if it is new. Since $f$ is unbounded below, we say that
a method to minimize $f$ is successful if it generates function values diverging to $-\infty$ and fails otherwise.

Theorem~3.1 of~\cite{AO2021} proves convergence of the absolute normalized search directions of the scaled method
when $a\ge2\sqrt r$. With some additional assumptions on the line search, this guarantees failure of
L-BFGS-1. The experiments in Section~4 of that paper suggest
that $a\ge\sqrt{3r}$ suffices. Section~\ref{sec:scaled} of this article proves this
stronger result and gives a geometric error estimation for $a \geq\sqrt{3r}$. 
These results also apply to the corresponding Theorem~3.2.1 in the first author's
dissertation~\cite{Asl2020}.

Section~\ref{sec:unscaled} addresses the contrasting behavior of the unscaled
method. Figure~5 of~\cite{AO2021} reports no failures for unscaled
memoryless BFGS in the experiments considered there, and Section~5 leaves
open whether the unscaled method can converge to a nonoptimal point.
We resolve the application of scaled L-BFGS-1 to~\eqref{eq:model} with an
Armijo--Wolfe line search; a bounded-below infinite run is impossible.
The result holds for all $a\ge\sqrt r$. 

These results form an addendum to~\cite{AO2021}; the original theorems
remain valid. The new proof replaces the monotone-subsequence argument
of  Appendix~B with a direct estimate of the error relative to a fixed
point. This particular convergence argument was not familiar to the authors
at the time. The two proofs arose through the
first author's directed discussions with GPT 6 agents. The proofs were later corrected, verified and rephrased by the authors.

\subsection{Scaled and unscaled memoryless BFGS}

Write $g_k=g(x_k)$ as defined in  \eqref{eq:gradient} and start with $x_0^{(1)}\ne0$ and $d_0=-g_0$.
After accepting a step $t_{k-1}>0$, define
\begin{equation}
 x_k=x_{k-1}+t_{k-1}d_{k-1},\qquad
 s_{k-1}=x_k-x_{k-1},\qquad y_{k-1}=g_k-g_{k-1}.
 \label{eq:pairs}
\end{equation}
For $k\ge1$, memoryless BFGS computes
\begin{align}
 \rho_{k-1}&=(s_{k-1}^Ty_{k-1})^{-1},\qquad
 V_{k-1}=I-\rho_{k-1}y_{k-1}s_{k-1}^T, \notag\\
 H_k&=V_{k-1}^T H_k^0 V_{k-1}
       +\rho_{k-1}s_{k-1}s_{k-1}^T,\qquad d_k=-H_kg_k .
 \label{eq:update}
\end{align}
Thus only the latest correction pair is used. The two choices considered
here are
\begin{equation}
 H_k^0=
 \begin{cases}
 \gamma_k I,\quad
 \displaystyle\gamma_k=\frac{s_{k-1}^Ty_{k-1}}{y_{k-1}^Ty_{k-1}},
       &\text{scaled},\\[6pt]
 I,    &\text{unscaled}.
 \end{cases}
 \label{eq:scaling}
\end{equation}
These are the standard definitions; see Nocedal and
Wright~\cite[Section~7.2, equations~(7.20) and~(7.22)]{NW2006}.
In particular, unscaled L-BFGS resets the initial matrix to the identity at every
iteration.

\subsection{Line search conditions and basic facts}
\label{sec:linesearch}

We use the convention of~\cite[Section~2]{AO2021}, that a step $t_k>0$
satisfies the Armijo-Wolfe conditions if
\begin{align}
 f(x_k+t_kd_k)&\le f(x_k)+c_1t_kg_k^Td_k,
       \label{eq:armijo}\\
 g(x_k+t_kd_k)&\text{ exists, and } 
 g(x_k+t_kd_k)^Td_k\ge c_2g_k^Td_k .
       \label{eq:wolfe}
\end{align}
with $0<c_1<c_2<1$.

As in~\cite[Section~2]{AO2021}, we say that
memoryless BFGS succeeds if it identifies the fact that $f$ is unbounded below  either because it generates a direction $d$ for
which no steplength $t$ satisfying the Armijo-Wolfe conditions exists (in which case the function must be unbounded below)  or, alternatively, that it generates a 
sequence $\{x_k\}$ for which Armijo-Wolfe steps always exist, with $f(x_k)\downarrow -\infty$.  If neither event takes place, $\{f(x_k)\}$ is bounded below,
which is regarded as failure, since $f$ is unbounded below.


\section{A sharpened convergence bound for scaled L-BFGS-1}
\label{sec:scaled}

For scaled L-BFGS-1, \cite[Lemma~2.1]{AO2021} shows that if $ a \ge \sqrt{3r}$, then 
\[
 b_k:=\frac{d_k^{(i)}}{d_k^{(1)}}\quad (i=2,\ldots,n)
\]
alternates in sign
with $|b_k|\leq 1/a$ and that $\beta_k:=|b_k|$ satisfies
\begin{equation}
 \beta_0=\frac1a,\qquad
 \beta_{k+1}=T(\beta_k),\qquad
 T(\beta):=\frac{1+r\beta^2}{a-r\beta}-\beta
      =\frac{1-a\beta+2r\beta^2}{a-r\beta}.
 \label{eq:scalar}
\end{equation}
Note that the recurrence in \eqref{eq:scalar} is exactly equation~(27) of~\cite{AO2021}.
Then \cite[Theorem~2.3]{AO2021} guarantees that consequently, provided  $ a \ge \sqrt{3r}$, the
L-BFGS-1 method generates a sequence of iterates and search directions for which
Armijo--Wolfe steps always exist and hence the method is well-defined.

Theorem~3.1 of~\cite{AO2021} proves that the sequence $\{\beta_k\}=\{|b_k|\}$ converges when $a\ge2\sqrt r$.
We now weaken this assumption to $a\ge\sqrt{3r}$ for any Armijo--Wolfe line search.
At the end of this section we state and justify the resulting failure guarantees, distinguishing arbitrary Armijo--Wolfe line searches from the specific bracketing search in~\cite[Algorithm~2]{AO2021}.

\begin{theorem}[Sharpening of {\cite[Theorem~3.1]{AO2021}}]
\label{thm:scaled}
Let $a\ge\sqrt{3r}$. The sequence $\beta_{k+1}=T(\beta_k)$ converges to
\begin{equation}
 b=\frac{a-\sqrt{a^2-3r}}{3r}
   =\frac1{a+\sqrt{a^2-3r}}.
 \label{eq:fixedpoint}
\end{equation}
\end{theorem}

\begin{proof}
As is shown in \cite[Lemma~2.1]{AO2021}, $0<T(\beta_k)\le1/a$ for $k=0,1,\hdots$, so all iterates lie in
 $J := [0,~1/a]$. Clearly $b$ in~\eqref{eq:fixedpoint} satisfies $0<b\le1/a$ and $T(b)=b$. That is, $b$ is a fixed-point of $T(\beta)$. It is easy to see that $1-ab=-3rb^2 + ab$ with some rearrangement of \eqref{eq:fixedpoint}. Substituting $1-ab$ 
by $-3rb^2 + ab$ as well as $rb\beta$ by $3rb\beta-2rb\beta$ in the numerator of
\begin{align*}
  T(\beta)-b =\frac{2r\beta^2+rb\beta -a\beta+1-ab}{a-r\beta}
\end{align*}
we arrive at the following factorization
\begin{align}
 T(\beta)-b
 &=(\beta-b)\,\frac{2r\beta+3rb-a}{a-r\beta}.
 \label{eq:factorization}
\end{align}
Now define the right factor in \eqref{eq:factorization} as 
$$R(\beta)=\dfrac{2r\beta+3rb-a}{a-r\beta}.$$
Hence, $T(\beta)-b = (\beta - b)R(\beta)$. For $\beta\in J$, it is easy to check that 
\[
1+R(\beta)=\frac{r(\beta+3b)}{a-r\beta}>0.
\]
Suppose first that $a>\sqrt{3r}$. Since $\beta,b\le1/a$,
\[
 3r(\beta+b)\le\frac{6r}{a}<2a,
\]
and hence
\[
 1-R(\beta)=\frac{2a-3r(\beta+b)}{a-r\beta}>0.
\]
Together with $1+R(\beta)>0$, this shows that $|R(\beta)|<1$ throughout $J$.
Moreover, $a-r\beta\ge a-r/a>0$ on $J$, so $R$ is continuous on the compact interval $J$.
Therefore its maximum absolute value is attained and satisfies
\[
 \kappa:=\max_{\beta\in J}|R(\beta)|<1.
\]
Applying~\eqref{eq:factorization} repeatedly gives
\[
 \abs{\beta_k-b}\le\kappa^k\abs{\beta_0-b}\longrightarrow0,
\]
which proves convergence for $a>\sqrt{3r}$.

If $a=\sqrt{3r}$, then $b=1/a$ and $R(b)=1$, so the preceding strict contraction bound does not apply.
However, the prescribed initial value in~\eqref{eq:scalar} is $\beta_0=1/a=b$.
Since $T(b)=b$, it follows that $\beta_k=b$ for every $k\ge0$.
Thus convergence also holds in this boundary case.  
\end{proof}

Theorem~\ref{thm:scaled} applies to \emph{any} line search satisfying~\eqref{eq:armijo}--\eqref{eq:wolfe}; it does not require the bracketing search of~\cite[Algorithm~2]{AO2021}.
Convergence of $\beta_k$ alone does not establish failure of the method.
The two distinct sufficient conditions for failure in~\cite[Theorems~3.8 and~3.11]{AO2021} extend to $a\ge\sqrt{3r}$ as follows.

For an arbitrary Armijo--Wolfe line search, add the additional condition
\begin{equation}
 \frac{1}{c_1}-1<\frac{a}{rb},
 \qquad\text{equivalently}\qquad c_1>\frac{rb}{a+rb}.
 \label{eq:scaled-failure-armijo}
\end{equation}
To see why this suffices, set $u_k:=\abs{s_k^{(1)}}$,
\[
 \varphi_k:=\frac{c_1(a+r\beta_k)+a-r\beta_k}{2a}
\quad\mathrm{and}\quad
 \varphi:=\frac{c_1(a+rb)+a-rb}{2a}.
\]
The Armijo inequality gives $\varphi_k u_k\le\abs{x_k^{(1)}}$, see~\cite[Lemma~3.3]{AO2021}.
Adding this inequality at two consecutive iterations and using the sign change
$u_{k-1}=\abs{x_{k-1}^{(1)}}+\abs{x_k^{(1)}}$ yields, as in~\cite[Corollary~3.5]{AO2021},
\[
 u_k\le\frac{1-\varphi_{k-1}}{\varphi_k}\,u_{k-1},
 \qquad
 \frac{1-\varphi_{k-1}}{\varphi_k}\longrightarrow\frac{1-\varphi}{\varphi}<1,
\]
where the final inequality follows from~\eqref{eq:scaled-failure-armijo}, since $\varphi-1/2=(c_1(a+rb)-rb) / (2a)>0$.
Thus there are $q\in(0,1)$ and $K$ such that $u_k\le q u_{k-1}$ for every $k\ge K$, and hence $\sum_{k=0}^{\infty}u_k<\infty$.
This proves the extension of~\cite[Theorem~3.8]{AO2021} to $a\geq \sqrt{3r}$, including $a=\sqrt{3r}$, when~\eqref{eq:scaled-failure-armijo} reduces to $c_1>1/4$.

For an Armijo--Wolfe search satisfying $t_k\le2$ for $k\ge1$, no additional restriction on $c_1$ is needed.
Indeed, the scaled update gives
\begin{align*}
 u_{k+1} &= t_{k+1}|d_{k=1}^{(1)}|= \frac{t_{k+1}}{2}\left(1-\frac{r\beta_k}{a}\right)u_k
       \le\left(1-\frac{r\beta_k}{a}\right)u_k,\\
 1-\frac{r\beta_k}{a}&\longrightarrow1-\frac{rb}{a}<1,
\end{align*}
where the second equality is the same as in~\cite[equation~(55)]{AO2021}.
Again we conclude $\sum_{k=0}^{\infty}u_k<\infty$.
The specific bracketing search in~\cite[Algorithm~2]{AO2021} has this steplength bound when $a\ge\sqrt{3r}$:
$\beta_k\le1/a$ implies $3r\beta_k\le a$, so the proofs of~\cite[Lemmas~3.4, 3.9 and~3.10]{AO2021} remain valid.
This proves the extension of~\cite[Theorem~3.11]{AO2021} to the sharpened range, independently of the Armijo parameter.

Under either sufficient condition, $\abs{s_k^{(i)}}=\beta_k u_k\le u_k/a$ for $i=2,\ldots,n$, so
\[
 \sum_{k=0}^{\infty}\norm{s_k}
 \le\sqrt{1+\frac{r}{a^2}}\sum_{k=0}^{\infty}u_k<\infty.
\]
Consequently $x_k$ converges to a finite point $\bar x$.
The sign changes give $\abs{x_k^{(1)}}\le u_k\to0$, so $\bar x^{(1)}=0$.
By continuity, $f(x_k)\downarrow f(\bar x)>-\infty$, establishing failure since $f$ is unbounded below.
The iterate-convergence conclusion of~\cite[Corollary~3.12]{AO2021} also holds under either set of assumptions.

\section{Success of unscaled L-BFGS-1}
\label{sec:unscaled}

\begin{theorem}
\label{thm:unscaled}
If we apply unscaled L-BFGS-1 to~\eqref{eq:model}, defined by~\eqref{eq:pairs}--\eqref{eq:update} with $H_k^0=I$ using any Armijo--Wolfe line search, and starting from any point $x_0$ with 
$x_0^{(1)}\ne0$, the method always succeeds.
\end{theorem}

\begin{proof}
As in~\cite[Section~2]{AO2021}, for this problem memoryless BFGS succeeds if it identifies the fact that $f$ is unbounded below either because it generates a direction $d$ for which no steplength $t$ satisfying the Armijo-Wolfe conditions exists (in which case the function must be unbounded below) or, alternatively, that it generates a sequence 
$\{x_k\}$ for which Armijo-Wolfe steps always exist and $f(x_k) \downarrow -\infty$.

As in~\cite{AO2021}, we assume without loss of generality that the
initial point has a positive first component, $x_0^{(1)}>0$. Furthermore we assume the method generates a sequence ${x_k}$ for which Armijo-Wolfe steps $t_k$ always exist. We then proceed to prove that the resulting function values diverge to $-\infty$.

For the function~\eqref{eq:model}, the Wolfe condition~\eqref{eq:wolfe} forces
every accepted step to change the sign of the first component of the update $x_k$,
\[
 \sgn(x_k^{(1)})=-\sgn(x_{k-1}^{(1)}),
\]
see~\cite[equations~(13) and~(14)]{AO2021}, whose derivation does not
depend on the choice of $H_k^0$. Hence, as long as the steps are
accepted, $x_k^{(1)}$ has the sign $(-1)^k$, and $x_k = x_{k-1}+s_{k-1}$ implies 
$(-1)^ks_{k-1}^{(1)}=\abs{s_{k-1}^{(1)}}>0$. Additionally 
\begin{equation}
 g_k=\left[\begin{array}{c}(-1)^ka\\ \one\end{array}\right], 
 \qquad
 y_{k-1}=\left[\begin{array}{c}2(-1)^ka\\0\end{array}\right],
  \qquad
  s_{k-1}=t_{k-1}d_{k-1},
 \label{eq:gradk}
\end{equation}
so
\begin{equation}
 s_{k-1}^Ty_{k-1}=2a\abs{s_{k-1}^{(1)}}.
 \label{eq:unscaledpair}
\end{equation}
All tail components of each direction, and hence of each step, are
equal: this holds for $d_0=-g_0$ and is preserved by~\eqref{eq:update}.
Moreover, $d_{k-1}^{(1)}\ne0$ for an accepted step. Hence, as in Section~\ref{sec:scaled} and
\cite[equation~(19)]{AO2021}, we can write
\begin{equation}
 b_{k-1}=\frac{d_{k-1}^{(i)}}{d_{k-1}^{(1)}}
        =\frac{s_{k-1}^{(i)}}{s_{k-1}^{(1)}},
 \quad i=2,\ldots,n,
 \qquad\text{so that}\qquad
 s_{k-1}=s_{k-1}^{(1)}\left[\begin{array}{c}1\\ b_{k-1}\one\end{array}\right].
 \label{eq:signedratio}
\end{equation}
For $k\ge1$, let $\gamma_k$ be the scaling parameter of~\eqref{eq:scaling},
which by~\eqref{eq:unscaledpair} is
\begin{equation}
 \gamma_k=\frac{s_{k-1}^Ty_{k-1}}{y_{k-1}^Ty_{k-1}}
         =\frac{\abs{s_{k-1}^{(1)}}}{2a}.
\label{eq:gamma}         
\end{equation}
Substituting this and~\eqref{eq:signedratio}
into~\eqref{eq:update} yields
\begin{equation}
 H_k=\left[
 \begin{array}{c|c}
  \gamma_k & \gamma_kb_{k-1}\one^T \\
  \hline
  \gamma_kb_{k-1}\one & (1+\gamma_k)b_{k-1}^2\one\one^T+I_r
 \end{array}
 \right],
 \label{eq:unscaledmatrix}
\end{equation}
where $I_r$ is the $r\times r$ identity matrix.
Although the unscaled method does not use $\gamma_k$ to scale $H_k^0$, the same quantity appears as a coefficient in the
updated matrix. From $d_k=-H_kg_k$, using ~\eqref{eq:gradk} we get 
\begin{equation}
 d_k=-\left[\begin{array}{c}
  \gamma_k\bigl((-1)^ka+rb_{k-1}\bigr)\\[4pt]
  \Bigl((-1)^ka\gamma_kb_{k-1}+(1+\gamma_k)rb_{k-1}^2+1\Bigr)\one
 \end{array}\right] = -d_k^{(1)}\left[\begin{array}{c}
 1\\
 b_k\one
 \end{array}\right].
 \label{eq:unscaleddir}
\end{equation}
In particular, note that the sign change implies that $d_k^{(1)}$ must have the same sign as $(-1)^{k+1}$, i.e.
\begin{equation}
 a+(-1)^krb_{k-1}>0.
 \label{eq:towardkink}
\end{equation}
If $t_{k+1}$ is accepted, using the same argument for iteration $k+1$ as we did for iteration $k$ we also get 
\begin{equation}
 a-(-1)^krb_k>0.
 \label{eq:towardkinkp1}
\end{equation}
Also from $s_k=t_kd_k$ we obtain 
\begin{equation}
 \abs{s_k^{(1)}}=t_k\gamma_k\bigl(a+(-1)^krb_{k-1}\bigr). \label{eq:srecurrence}
\end{equation}
Recall   
\begin{align}
 b_k &=\frac{(-1)^ka\gamma_kb_{k-1}+(1+\gamma_k)rb_{k-1}^2+1}{\gamma_k\bigl((-1)^ka+rb_{k-1}\bigr)} \notag \\
 &=b_{k-1}+(-1)^k\,\frac{1+rb_{k-1}^2}{\gamma_k\bigl(a+(-1)^krb_{k-1}\bigr)}.  \label{eq:brecurrence}
\end{align}
Plugging $b_k$ into~\eqref{eq:towardkinkp1} we obtain
\begin{equation*}
 a - (-1)^krb_{k-1} > \frac{r(1+{r}b_{k-1}^2)}{\gamma_k\bigl(a+(-1)^krb_{k-1}\bigr)},
\end{equation*}
The right term is positive by~\eqref{eq:towardkink}, so this inequality shows
\begin{equation}
 a+(-1)^krb_{k-1}<2a,
 \label{eq:kinkupper}
\end{equation}
so 
\begin{equation*}
 1 > \frac{r(1+{r}b_{k-1}^2)}{\gamma_k\bigl( a - (-1)^krb_{k-1} \bigr) \bigl(a+(-1)^krb_{k-1}\bigr)}.
\end{equation*}
Moving $\gamma_k$ to the left side and using \eqref{eq:gamma} we get
\begin{equation*}
\abs{s_{k-1}^{(1)}} > \dfrac{2ar(1+rb_{k-1}^2)}{\Bigl(a+(-1)^krb_{k-1}\Bigr)\Bigl(a-(-1)^krb_{k-1}\Bigr)}.
\end{equation*}
Since $s_{k-1}$ is in the direction of descent, $-t_{k-1}g_{k-1}^Td_{k-1} = -g_{k-1}^Ts_{k-1} = \abs{g_{k-1}^Ts_{k-1}}$ and 
\begin{align*}
\abs{g_{k-1}^Ts_{k-1}} &= (-1)^{{k}}s_{k-1}^{(1)}(a + (-1)^{k-1}rb_{k-1})\\
&= \abs{s_{k-1}^{(1)}}(a + (-1)^{k-1}rb_{k-1})\\
&> \dfrac{2ar(1+rb_{k-1}^2)\Bigl(a + (-1)^{k-1}rb_{k-1}\Bigr)}{\Bigl(a-(-1)^{k-1}rb_{k-1}\Bigr)\Bigl(a+(-1)^{k-1}rb_{k-1}\Bigr)}\\
&= \dfrac{2ar(1+rb_{k-1}^2)}{\Bigl(a +(-1)^{k}rb_{k-1}\Bigr)} \\
&> \dfrac{2ar(1+rb_{k-1}^2)}{2a} \ge r.
\end{align*}
The second equality comes from recalling $\abs{s_{k-1}^{(1)}}=(-1)^ks_{k-1}^{(1)}$. The last strict inequality is from~\eqref{eq:kinkupper}.

By the Armijo condition~\eqref{eq:armijo},
\[
 f(x_{k+1})\le f(x_k)+c_1t_kg_k^Td_k
 \le f(x_k)-c_1r,
\]
where the second inequality uses $-t_kg_k^Td_k\ge r$.
Summing over $k=1,\ldots,N-1$ gives
\[
 f(x_N)\le f(x_1)-Nc_1r\longrightarrow-\infty,
\]
since $c_1>0$ and $r=n-1\ge1$.

Thus, if Armijo--Wolfe steps exist at every iteration, the function values diverge to $-\infty$.
If instead no Armijo--Wolfe step exists at some iteration, the method succeeds by identifying a direction along which $f$ is unbounded below, as discussed in Section~\ref{sec:linesearch}.
These two alternatives establish that the unscaled method always succeeds.

\end{proof}


\begin{thebibliography}{9}

\bibitem{Asl2020}
A.~Asl.
\newblock \emph{Behavior of the Limited-Memory BFGS Method on Nonsmooth
Optimization Problems in Theory and Practice}.
\newblock Ph.D. dissertation, Department of Computer Science,
New York University, May 2020.
\newblock \url{https://cs.nyu.edu/media/publications/asl_thesis_final_UtpoLsu.pdf}.

\bibitem{AO2021}
A.~Asl and M.~L. Overton.
\newblock Analysis of limited-memory BFGS on a class of nonsmooth convex functions.
\newblock \emph{IMA Journal of Numerical Analysis}, 41(1):1--27, 2021.
\newblock \href{https://doi.org/10.1093/imanum/drz052}{doi:10.1093/imanum/drz052}.

\bibitem{NW2006}
J.~Nocedal and S.~J. Wright.
\newblock \emph{Numerical Optimization}.
\newblock Springer Series in Operations Research and Financial Engineering.
\newblock Springer, New York, second edition, 2006.
\newblock \href{https://doi.org/10.1007/978-0-387-40065-5}{doi:10.1007/978-0-387-40065-5}.

\end{thebibliography}
\end{document}